\documentclass{aims}
\usepackage{amsmath}
  \usepackage{paralist}
  \usepackage{graphics}
  \usepackage{epsfig}
\usepackage{graphicx}  \usepackage{epstopdf}
 \usepackage[colorlinks=true]{hyperref}
\hypersetup{urlcolor=blue, citecolor=red}

\usepackage{url}

\makeatletter
\def\publname{}
\def\volinfo{}
\def\pageinfo{}
\makeatother

\usepackage{xspace}
\usepackage[shortlabels]{enumitem}
\usepackage{dsfont}
\usepackage{amsopn}

\newcommand\indic{\mathds{1}}
\newcommand\aew{\mathrm{ a.e. } \,\, }

\newcommand{\lin}[1]{\mathcal{L}(#1)}
\def\norm#1{\left\|#1\right\|}

\newcommand{\ens}[1]{ \left\{#1\right\} }
\newcommand{\dom}[1]{ D(#1) }
\newcommand{\field}[1]{\ensuremath{\mathbb{#1}}}
\newcommand{\C}{\field{C}\xspace}
\newcommand{\R}{\field{R}\xspace}

\newcommand{\N}{\field{N}\xspace}
\newcommand\ds{\displaystyle}
\newcommand\Id{\mathrm{Id}}
\newcommand\ddt{\frac{d}{dt}}

\newcommand\opA{A}
\newcommand\opB{B}

\newcommand{\K}{K}
\newcommand\Tmin{T_0}
\newcommand\range{\mathrm{Im} \,}
\newcommand\spec[1]{\sigma\left(#1\right)}

\newtheorem{theorem}{Theorem}[section]
\newtheorem{corollary}[theorem]{Corollary}

\newtheorem{lemma}[theorem]{Lemma}

\theoremstyle{definition}

\newtheorem{remark}[theorem]{Remark}

\title[Perturbations of controlled systems]
      {Compact perturbations of controlled systems}

\author[Michel Duprez and Guillaume Olive]{}


 \email{mduprez@math.cnrs.fr}
 \email{math.golive@gmail.com}

\thanks{$^*$ Corresponding author: Guillaume Olive}

\begin{document}
\maketitle

\centerline{\scshape Michel Duprez}
\medskip
{\footnotesize
 \centerline{ Institut de Math\'ematique de Marseille}
 \centerline{Aix Marseille Universit\'e}
 \centerline{39, rue J.Joliot Curie, 13453 Marseille Cedex 13, France}
}

\medskip

\centerline{\scshape Guillaume Olive$^*$}
\medskip
{\footnotesize
 \centerline{Institut de Math\'ematiques de Bordeaux}
 \centerline{Universit\'e de Bordeaux}
 \centerline{351, Cours de la Lib\'eration, 33405 Talence, France}
}

\bigskip

\begin{abstract}
In this article we study the controllability properties of general compactly perturbed exactly controlled linear systems with admissible control operators.
Firstly, we show that approximate and exact controllability are equivalent properties for such systems.
Then, and more importantly, we provide for the perturbed system a complete characterization of the set of reachable states in terms of the Fattorini-Hautus test.
The results rely on the Peetre lemma.
\end{abstract}

\section{Introduction and main results}

In this work we study the exact controllability property of general compactly perturbed controlled linear systems using a compactness-uniqueness approach.
This technique was notably used in the pioneering work \cite{RT74} to establish the exponential decay of the solutions to some hyperbolic equations.
On the other hand, the first controllability results using this method were obtained in \cite{Z87,L88} for a plate equation and in \cite{Z91} for a perturbed wave equation.
Wether one wants to establish a stability result or a controllability result, one is led in both cases to prove estimates, energy estimates or observability inequalities.
For a perturbed system, a general procedure is to start by the known estimate satisfied by the unperturbed system and to try to derive the desired estimate, up to some ``lower order terms'' that we would like to remove.
The compactness-uniqueness argument then reduces the task of absorbing these additional terms to a unique continuation property for the perturbed system.
We should point out that, despite the numerous applications of this flexible method to successfully establish the controllability of systems governed by partial differential equations (see e.g. \cite{Z87,L88,Z91,BLR92,LT05,CT10,FCLZ16,LR17}, etc.), no systematic treatment has been provided so far for compactly perturbed systems, by which we mean that there is no abstract result available in the literature that covers all type of systems, regardless the nature of the equations we are considering (plate, wave, etc.).
This will be the first point of the present paper to fill this gap (Theorem \ref{thm pert AC} below).
Then, and more importantly, we considerably improve this result by establishing an explicit characterization of the set of reachable states for the perturbed system (Theorem \ref{thm pert 2} below).
This characterization is given in terms of the Fattorini-Hautus test, which is a far weaker kind of unique continuation property than the approximate controllability.
Our result shows in particular that this test is actually sufficient to ensure the exact controllability of the perturbed system (Corollary \ref{thm pert fatt} below).
We illustrate this corollary and the significance of the Fattorini-Hautus test with an application to the exact controllability of a Schr\"odinger equation perturbed by a non local term (Theorem \ref{thm schro} below).
Finally, we conclude this work by mentioning that the techniques used in this paper are not only limited to the study of the controllability of compactly perturbed controlled systems (Theorem \ref{main thm} below).
The proofs of the main results of this article are based on the Peetre Lemma, introduced in \cite{P61}, which is in fact the root of compactness-uniqueness methods.

Let us now introduce some notations and recall some basic facts about the controllability of abstract linear evolution equations.
We refer to the excellent textbook \cite{TW09} for the proofs of the statements below.
Let $H$ and $U$ be two (real or complex) Hilbert spaces, let $\opA:\dom{\opA} \subset H \longrightarrow H$ be the generator of a $C_0$-semigroup $(S_{\opA}(t))_{t \geq 0}$ on $H$ and let $\opB \in \lin{U,\dom{\opA^*}'}$.
For $T \geq 0$ let $\Phi_T \in \lin{L^2(0,+\infty;U),\dom{\opA^*}'}$ be the input map of $(A,B)$, that is the linear operator defined for every $u \in L^2(0,+\infty;U)$ by
$$\Phi_T u=\int_0^T S_{A}(T-s)Bu(s) \, ds.$$
We assume that $\opB$ is admissible for $\opA$, which means that $\range \Phi_T \subset H$ for some (and hence all) $T>0$.
From this assumption it follows that $\Phi_T \in \lin{L^2(0,+\infty;U),H}$.
Its adjoint $\Phi_T^* \in \lin{H,L^2(0,+\infty;U)}$ is the unique continuous linear extension to $H$ of the map $z \in \dom{A^*} \mapsto B^*S_{\opA}(T-\cdot)^*z \in L^2(0,+\infty;U)$, where $B^*S_{\opA}(T-\cdot)^*z$ is extended by zero outside $(0,T)$ (in particular, $\Phi_T^*z(t)=0$ for a.e. $t>T$).
Let us now consider the abstract evolution system
\begin{equation}\label{syst abst}
\left\{
\begin{array}{rcll}
\ds \ddt y&=& \opA y + \opB u, & t \in (0,T), \\
y(0) &=& y^0, &
\end{array}
\right.
\end{equation}
where $T>0$ is the time of control, $y^0 \in H$ is the initial data, $y$ is the state and $u \in L^2(0,T;U)$ is the control.
Since $\opB$ is admissible for $\opA$, system \eqref{syst abst} is well-posed: for every $y^0 \in H$ and every $u \in L^2(0,T;U)$, there exists a unique solution $y \in C^0([0,T];H)$ to system \eqref{syst abst} given by the Duhamel formula
$$y(t)=S_A(t)y^0+\Phi_t u, \quad \forall t \geq 0.$$
The regularity of the solution allows us to consider control problems for system \eqref{syst abst}.
We say that system \eqref{syst abst} or $(\opA,\opB)$ is:
\begin{itemize}
\item
exactly controllable in time $T$ if, for every $y^0,y^1 \in H$, there exists $u \in L^2(0,T;U)$ such that the corresponding solution $y$ to system \eqref{syst abst} satisfies $y(T)=y^1$.
\item
approximately controllable in time $T$ if, for every $\varepsilon>0$ and every $y^0,y^1 \in H$, there exists $u \in L^2(0,T;U)$ such that the corresponding solution $y$ to system \eqref{syst abst} satisfies $\norm{y(T)-y^1}_{H} \leq \varepsilon$.
\end{itemize}
Clearly, exact controllability in time $T$ implies approximate controllability in the same time.
The set $\range \Phi_T$ (resp. $\overline{\range \Phi_T}$) is called the set of exactly (resp. approximately) reachable states in time $T$.
Therefore, $(A,B)$ is exactly (resp. approximately) controllable in time $T$ if, and only if, $\range \Phi_T=H$ (resp. $\overline{\range \Phi_T}=H$).
It is also well-known that the controllability has a dual concept named observability.
More precisely, $(\opA,\opB)$ is exactly controllable in time $T$ if, and only if, there exists $C>0$ such that
\begin{equation}\label{obs}
\norm{z}_H^2 \leq C \int_0^T \norm{\Phi_T^*z(t)}_U^2 \, dt, \quad \forall z \in H,
\end{equation}
and $(\opA,\opB)$ is approximately controllable in time $T$ if, and only if,
\begin{equation}\label{ucp}
\Big(\Phi_T^*z(t)=0, \quad \aew t \in (0,T)\Big) \Longrightarrow z=0, \quad \forall z \in H.
\end{equation}

Let us now state the main results of this paper.
The first one simply unifies previous results available in the literature under a general semigroup setting:

\begin{theorem}\label{thm pert AC}
Let $H$ and $U$ be two (real or complex) Hilbert spaces.
Let $\opA_0:\dom{\opA_0} \subset H \longrightarrow H$ be the generator of a $C_0$-semigroup on $H$ and let us consider $\opB \in \lin{U,\dom{\opA_0^*}'}$ an admissible control operator for $\opA_0$.
Let $\K \in \lin{H}$ and let us form the unbounded operator $\opA_{\K}=\opA_0+\K$ with $\dom{\opA_{\K}}=\dom{\opA_0}$
\footnote{
$\opA_{\K}$ is then the generator of a $C_0$-semigroup on $H$ and $\opB$ is also admissible for $\opA_{\K}$, see below.
}.
We assume that:
\begin{enumerate}[(i)]
\item\label{hyp1 pAC}
There exists $\Tmin>0$ such that $(\opA_0,\opB)$ is exactly controllable in time $\Tmin$.
\item\label{hyp2 pAC}
$\K$ is compact.
\item\label{item ac}
$(\opA_{\K},\opB)$ is approximately controllable in time $\Tmin$.
\end{enumerate}
Then, $(\opA_{\K},\opB)$  is exactly controllable in time $\Tmin$.
\end{theorem}

The second and most important result of the present paper shows that we can even give a very precise characterization of the reachable states for the perturbed system, if we allow the time of control to be slightly longer:

\begin{theorem}\label{thm pert 2}
Let $H$ and $U$ be complex Hilbert spaces and let $A_0, A_K$ and $B$ be defined as in Theorem \ref{thm pert AC}.
For $T \geq 0$ let $\Phi_T \in \lin{L^2(0,+\infty;U),H}$ be the input map of $(A_K,B)$.
Let $\sigma_F$ be the set given by
$$\sigma_F=\ens{\lambda \in \C, \quad \ker(\lambda-\opA_K^*) \cap \ker \opB^* \neq \ens{0}},$$
and for every $\lambda \in \C$ let $E_\lambda$ be the subspace of $H$ defined by
$$E_\lambda=\ens{z \in \bigcup_{m=1}^{+\infty} \ker(\lambda-A_K^*)^m, \quad B^*(\lambda-A_K^*)^m z=0, \quad \forall m \in \N}.$$
Then, under the assumptions \ref{hyp1 pAC} and \ref{hyp2 pAC} of Theorem \ref{thm pert AC}, the set $\sigma_F$ is finite, $E_\lambda$ is finite dimensional for every $\lambda \in \C$, and we have
$$\range \Phi_T =\left(\bigoplus_{\lambda \in \sigma_F} E_\lambda\right)^{\perp}, \quad \forall T>\Tmin.$$
\end{theorem}

This second result shows in particular that the approximate controllability assumption \ref{item ac} of Theorem \ref{thm pert AC} can be weakened to the Fattorini-Hautus test:

\begin{corollary}\label{thm pert fatt}
Let $H$ and $U$ be complex Hilbert spaces and let $A_0, A_K$ and $B$ be defined as in Theorem \ref{thm pert AC}.
Then, under the assumptions \ref{hyp1 pAC} and \ref{hyp2 pAC} of Theorem \ref{thm pert AC}, and if $(\opA_{\K},\opB)$ satisfies the Fattorini-Hautus test:
\begin{equation}\label{fatt test}
\ker(\lambda-\opA_{\K}^*) \cap \ker \opB^*=\ens{0}, \quad \forall \lambda \in \C,
\end{equation}
then $(\opA_{\K},\opB)$  is exactly controllable in time $T$ for every $T>\Tmin$.
\end{corollary}

\begin{remark}
The assumption \ref{hyp2 pAC} can be relaxed in some cases, see Remark \ref{rem lem cpct} below.
Besides, the abstract method we will develop in this paper is not only limited to the study of the controllability of perturbed controlled systems and it can be used to obtain the exact same conclusions as in the previous theorems for some (unperturbed) systems themselves.
We explained it in Section \ref{sect 4} below.
\end{remark}

\begin{remark}
It follows from Corollary \ref{thm pert fatt} that, if $(A_0,B)$ and $(A_K,B)$ are two systems satisfying the Fattorini-Hautus test, and $K$ is compact, then
\begin{multline*}
\inf\ens{T>0, \,\, (A_0,B) \mbox{ is exactly controllable in time } T}
\\ =\inf\ens{T>0, \,\, (A_K,B) \mbox{ is exactly controllable in time } T}.
\end{multline*}
In other words, both systems share the same minimal time of control.
\end{remark}

Corollary \ref{thm pert fatt} shows that, in order to prove the exact controllability of a compactly perturbed system which is known to be exactly controllable, it is (necessary and) sufficient to only check the Fattorini-Hautus test \eqref{fatt test}.
This result has been established in a particular case in \cite[Theorem 5]{CT10} for a perturbed plate equation with distributed controls.
The Fattorini-Hautus test appears for the very first time in \cite[Corollary 3.3]{F66} and it is also sometimes misleadingly known as the Hautus test in finite dimension, despite it has been introduced earlier by Fattorini, moreover in a much larger setting.
In a complete abstract control theory framework, it is the sharpest sufficient condition one can hope for since it is always a necessary condition for the exact, null or approximate controllability, to hold in some time (this easily follows from the dual characterizations \eqref{obs} or \eqref{ucp}).
It is also nowadays well-known that this condition characterizes the approximate controllability of a large class of systems generated by analytic semigroups (see \cite{F66,BT14}).
In practice, the Fattorini-Hautus test can be checked by various techniques, such as Carleman estimates for stationary systems (see e.g. \cite{CT10,BT14}), or spectral analysis when this later technique is not available (see e.g. \cite{BO14,CHO16} or the example of Section \ref{sect 3} below).

Let us mention that wether the Fattorini-Hautus test \eqref{fatt test} remains sufficient or not to obtain the exact controllability of the perturbed system in time $\Tmin$ may depend on the systems under consideration.
Therefore, both Theorem \ref{thm pert AC} and Corollary \ref{thm pert fatt} are important.
Obviously, Corollary \ref{thm pert fatt} is a stronger result if we do not look for the best time.
However, in some situations it might be important to preserve the same time of control (see e.g. \cite[Theorem 1.1]{CHO16}).

Finally, let us point out that throughout this work we do not request any spectral properties whatsoever on the operators $\opA_0$ or $\opA_{\K}$, contrary to the papers \cite{KL00,M04} where the existence of a Riesz basis of generalized eigenvectors or related spectral properties are required.

The rest of this paper is organized as follows.
In Section \ref{sect 2}, we prove the main results of this work.
In Section \ref{sect 3}, we illustrate the abstract results with an application to the controllability of a perturbed Schr\"odinger equation.
In Section \ref{sect 4}, we point out that the techniques used in this paper are not only restricted to the study of the controllability of compactly perturbed controlled systems.
Finally, we have included in Appendix \ref{appendix} a proof of an estimate that is needed in the proof of our main results (especially for unbounded admissible control operators).

\section{Proofs of the results}\label{sect 2}

The proofs of Theorem \ref{thm pert AC} and Theorem \ref{thm pert 2} both rely on the Peetre Lemma (see \cite[Lemma 3]{P61}):
\begin{lemma}\label{lemme peetre}
Let $B_1,B_2$ be two Banach spaces and let $L \in \lin{B_1,B_2}$.
The following two conditions are equivalent:
\begin{enumerate}[(a)]
\item\label{hyp a}
There exist a Banach space $B_3$, a compact operator $P \in \lin{B_1,B_3}$ and $\alpha>0$ such that
\begin{equation}\label{estim peetre}
\alpha\norm{z}_{B_1} \leq \norm{Lz}_{B_2}+\norm{Pz}_{B_3}, \quad \forall z \in B_1.
\end{equation}
\item\label{hyp b}
$\ker L$ is finite dimensional and $\range L$ is closed.
\end{enumerate}
\end{lemma}

\begin{remark}\label{rem lemme peetre}
It is well-known that $\ker L=\ens{0}$ and $\range L$ is closed if, and only if, there exists $\beta>0$ such that (see e.g. \cite[Lemma 4]{P61})
$$\beta\norm{z}_{B_1} \leq \norm{Lz}_{B_2}, \quad \forall z \in B_1.$$
Therefore, the compact term in \eqref{estim peetre} can be cancelled if one shows that $\ker L=\ens{0}$.
\end{remark}

Let us denote by $(S_{\opA_0}(t))_{t \geq 0}$ (resp. $(S_{\opA_{\K}}(t))_{t \geq 0}$) the $C_0$-semigroup generated by $\opA_0$ (resp. $\opA_{\K}$).
For $T \geq 0$ let $\Psi_T \in \lin{L^2(0,+\infty;U),H}$ (resp. $\Phi_T \in \lin{L^2(0,+\infty;U),H}$) be the input map of $(A_0,B)$ (resp. $(A_K,B)$).
Assume now that $(\opA_0,\opB)$ is exactly controllable in time $\Tmin$.
Then, for every $T \geq \Tmin$, there exists $C>0$ such that, for every $z \in H$,
$$\norm{z}_H^2 \leq C \int_0^T \norm{\Psi_T^*z(t)}_U^2 \, dt,$$
so that
\begin{equation}\label{obs estim cpct err}
\norm{z}_H^2 \leq 2C \left(\int_0^T \norm{\Phi_T^*z(t)}_U^2 \, dt+\int_0^T \norm{\Psi_T^*z(t)-\Phi_T^*z(t)}_U^2 \, dt\right).
\end{equation}
To prove that $(\opA_K,\opB)$ is exactly controllable in time $T$, we would like to get rid of the last term in the right-hand side of \eqref{obs estim cpct err}.
Therefore, we would like to apply Lemma \ref{lemme peetre} (see also Remark \ref{rem lemme peetre}) to the operators $L=\Phi_T^*$ and $P=\Psi_T^*-\Phi_T^*$.
Note that both operators are bounded linear operators since $\opB$ is admissible for both $\opA_0$ and $\opA_{\K}$ (see below).
To apply Lemma \ref{lemme peetre}, we have to check that $\Psi_T^*-\Phi_T^*$ is compact.

\begin{lemma}\label{lem cpct}
The operator $\Psi_T^*-\Phi_T^* \in \lin{H,L^2(0,+\infty;U)}$ is compact for every $T>0$.
\end{lemma}

For the proof of Lemma \ref{lem cpct} we need to recall the following estimate (see Appendix \ref{appendix}): for every $T>0$, there exists $C>0$ such that, for every $f \in C^1([0,T];H)$,
\begin{equation}\label{key estim}
\int_0^T \norm{\opB^* \int_0^t S_{\opA_0}(t-s)^*f(s) \,ds}_{U}^2 \, dt
\leq C \norm{f}_{L^2(0,T;H)}^2.
\end{equation}
This estimate holds because $B$ is admissible for $A_0$ (for bounded operators $B \in \lin{U,H}$ it is a straightforward consequence of the Cauchy-Schwarz inequality).
Using the dual characterization of admissibility (see \eqref{B admiss} below) and combining \eqref{key estim} with the identity \eqref{integr rela} below we also see that, if $B$ is admissible for $A_0$, then $B$ is admissible for $A_K$ as well.

\begin{proof}[Proof of Lemma \ref{lem cpct}]
Let us first compute $\Psi_T^*-\Phi_T^*$.
To this end, we recall the integral equation satisfied by semigroups of boundedly perturbed operators (see e.g. \cite[Corollary III.1.7]{EN}), valid for every $z \in H$ and $t \in [0,T]$:
\begin{equation}\label{integr rela}
S_{\opA_{\K}}(t)^*z=
S_{\opA_0}(t)^*z
+\int_0^t S_{\opA_0}(t-s)^*Fz(s) \,ds,
\end{equation}
where we introduced $F \in \lin{H,L^2(0,T;H)}$ defined for every $z \in H$ and every $s \in [0,T]$ by
$$Fz(s)= \K^*S_{\opA_{\K}}(s)^*z.$$
Note that $Fz \in C^1([0,T];H)$ for $z \in \dom{\opA_{\K}^*}$.
Thus, we have $\int_0^t S_{\opA_0}(t-s)^* Fz(s) \,ds \in \dom{A_0^*}$ for every $t \in (0,T)$ if $z \in \dom{\opA_{\K}^*}$.
This shows that each term in \eqref{integr rela} actually belongs to $\dom{\opA_{0}^*}$ if $z \in\dom{\opA_{0}^*}=\dom{\opA_{\K}^*}$.
Therefore, we can apply $\opB^*$ to obtain the following expression for $\Psi_T^*-\Phi_T^*$:
$$(\Psi_T^*-\Phi_T^*) z(t)=-\opB^* \int_0^t S_{\opA_0}(t-s)^* Fz(s) \,ds,$$
for every $z \in \dom{\opA_0^*}$ and a.e. $t \in (0,T)$.
Using now \eqref{key estim} there exists $C>0$ such that
$$\norm{(\Psi_T^*-\Phi_T^*) z}_{L^2(0,T;U)} \leq C \norm{Fz}_{L^2(0,T;H)},$$
for every $z \in \dom{\opA_{\K}^*}$, and thus for every $z \in H$ by density.
To conclude the proof it only remains to show that $F$ is compact.
Since $H$ is a Hilbert space, we will prove that, if $(z_n)_n \subset H$ is such that $z_n \to 0$ weakly in $H$ as $n \to +\infty$, then $F z_n \to 0$ strongly in $L^2(0,T;H)$ as $n \to +\infty$.
Since $z_n \to 0$ weakly in $H$ as $n \to +\infty$, using the strong (and therefore weak) continuity of semigroups on $H$, we obtain
$$S_{\opA_{\K}}(s)^*z_n \xrightarrow[n \to +\infty]{} 0 \quad \mbox{ weakly in } H, \quad \forall s \in [0,T].$$
Since $\K^*$ is compact, we obtain
$$\K^*S_{\opA_{\K}}(s)^*z_n \xrightarrow[n \to +\infty]{} 0 \quad \mbox{ strongly in } H, \quad \forall s \in [0,T].$$
On the other hand, by the classical semigroup estimate, $(\K^*S_{\opA_{\K}}(s)^*z_n)_n$ is clearly uniformly bounded in $H$ with respect to $s$ and $n$.
Therefore, the Lebesgue's dominated convergence theorem applies, so that $F z_n \to 0$ strongly in $L^2(0,T;H)$ as $n \to +\infty$.
This shows that $F$ is compact.
\end{proof}

\begin{remark}\label{rem lem cpct}
In all this work, the assumption \ref{hyp2 pAC} that $K$ is compact is only used to establish Lemma \ref{lem cpct}.
Therefore, the results obtained in Theorem \ref{thm pert AC} and \ref{thm pert 2} remain valid if, instead of assuming that $K$ is compact, we assume that the difference of the input maps is compact.
This assumption is for instance satisfied by some first order hyperbolic systems (see e.g. \cite[Thereom A]{NRL86}).
This assumption is also satisfied if $A_0$ generates an immediately compact semigroup, as it is easily seen from the previous proof by changing the roles of $A_0$ and $A_K$ and switching the convergence arguments at the end.
This concerns for instance the Korteweg-de Vries equation studied in \cite[Section 3]{R97}.
Let us point out that this is only interesting if we consider unbounded control operators since the exact controllability is impossible if the semigroup is immediately compact and the control operator is bounded (see \cite[Theorem 1.2]{T77}).
\end{remark}

The proof of Theorem \ref{thm pert AC} is now a direct consequence of Lemma \ref{lemme peetre} and \ref{lem cpct}.

\begin{proof}[Proof of Theorem \ref{thm pert AC}]
The assumptions of Lemma \ref{lemme peetre} are satisfied for $L=\Phi_{\Tmin}^*$ and $P=\Psi_{\Tmin}^*-\Phi_{\Tmin}^*$.
Therefore, $\range \Phi_{\Tmin}$ is closed (we recall that $\range \Phi_{\Tmin}$ is closed if, and only if, so is $\range \Phi_{\Tmin}^*$) and it follows from the very definitions of the notions of controllability that $(A_K,B)$ is then exactly controllable in time $\Tmin$ if, and only if, $(A_K,B)$ is approximately controllable in time $\Tmin$.
\end{proof}

\begin{remark}
We see that we have actually obtained a more precise conclusion than the one stated in Theorem \ref{thm pert AC}, namely that the targets that can be reached approximately can also be reached exactly (wether the whole system is controllable or not).
\end{remark}

Note that so far we have used only the second part of the conclusion of Lemma \ref{lemme peetre}.
For the proof of Theorem \ref{thm pert 2} we need the following general result:

\begin{lemma}\label{lem fatt}
Let $H$ and $U$ be two complex Hilbert spaces.
Let $\opA:\dom{\opA} \subset H \longrightarrow H$ be the generator of a $C_0$-semigroup on $H$ and let $\opB \in \lin{U,\dom{\opA^*}'}$ be admissible for $\opA$.
For $T \geq 0$ let $\Phi_T \in \lin{L^2(0,+\infty;U),H}$ be the input map of $(A,B)$.
Let $\sigma_F$ be the set given by
\begin{equation}\label{fatt spec}
\sigma_F=\ens{\lambda \in \C, \quad \ker(\lambda-\opA^*) \cap \ker \opB^* \neq \ens{0}},
\end{equation}
and for every $\lambda \in \C$ let $E_\lambda$ be the subspace of $H$ defined by
\begin{equation}\label{Elambda}
E_\lambda=\ens{z \in \bigcup_{m=1}^{+\infty} \ker(\lambda-A^*)^m, \quad B^*(\lambda-A^*)^m z=0, \quad \forall m \in \N}.
\end{equation}
Assume that there exists $\Tmin>0$ such that
\begin{equation}\label{dim ker PhiTmin finie}
\dim \ker \Phi_{\Tmin}^* <+\infty.
\end{equation}
Then, the set $\sigma_F$ is finite, $E_\lambda$ is finite dimensional for every $\lambda \in \C$, and we have
\begin{equation}\label{ker PhiT}
\ker \Phi_T^* =\bigoplus_{\lambda \in \sigma_F} E_\lambda, \quad \forall T>\Tmin.
\end{equation}
\end{lemma}

This result already appeared in the literature for some particular problems and in the particular case $\sigma_F=\emptyset$ (see e.g. \cite[Theorem 3.8]{BLR92}, \cite[Lemma 3.4]{R97}, \cite[Theorem 5]{CT10}, etc.).
To the best of our knowledge, the original ideas seem to trace back to \cite[Proposition 2]{RT74}.
However, we will follow the more recent presentation done in the proof of \cite[Theorem 5]{CT10}.

\begin{remark}
In the finite dimensional case $H=\C^n$ and $U=\C^m$ ($n,m \in \N^*$) we recover the well-known fact that
$\range \Phi_T=\range(B|AB|\cdots|A^{n-1}B)$ for every $T>0$.
\end{remark}

\begin{proof}[Proof of Lemma \ref{lem fatt}]
Let us first prove that, for every $T>0$ and $\lambda \in \C$, we have
\begin{equation}\label{Elambda inclus}
\ker \Phi_T^*  \supset E_\lambda.
\end{equation}
Let then $z \in E_\lambda$.
Thus, $z \in \dom{(A^*)^{\infty}}$ and there exists $m \in \N^*$ such that
\begin{equation}\label{first cond}
(\lambda-A^*)^m z=0,
\end{equation}
and
\begin{equation}\label{second cond}
B^*(\lambda-A^*)^r z=0, \quad \forall r \in \ens{0,\ldots,m-1}.
\end{equation}
Thanks to \eqref{first cond} we have, for every $t \geq 0$,
$$S_A(t)^*z=e^{\lambda t}\sum_{r=0}^{m-1} \frac{t^r}{r!} (A^*-\lambda)^r z.$$
Applying $B^*$ and using \eqref{second cond} we obtain that $z \in \ker \Phi_T^*$.
This establishes \eqref{Elambda inclus}.
Since the sum $\sum_{\lambda \in \sigma_F} E_\lambda$ is a direct sum, \eqref{Elambda inclus} implies that
\begin{equation}\label{ker PhiT supset}
\ker \Phi_T^*  \supset \bigoplus_{\lambda \in \sigma_F} E_\lambda, \quad \forall T>0.
\end{equation}
In particular, by \eqref{dim ker PhiTmin finie} we obtain that $\sigma_F$ is finite and that $E_\lambda$ is finite dimensional for every $\lambda \in \sigma_F$.

Let us now prove the reverse inclusion for $T>\Tmin$.
The key point is to establish that
\begin{equation}\label{in the dom}
\ker \Phi_T^* \subset \dom{\opA^*}.
\end{equation}
Let us first recall that, for every $T \geq T'$, we have $\ker \Phi_{T}^* \subset \ker \Phi_{T'}^*$ since
$$\int_0^{T'} \norm{\Phi_{T'}^*z(t)}_{U}^2 \, dt \leq \int_0^{T} \norm{\Phi_{T}^*z(t)}_{U}^2 \, dt, \quad \forall z \in H,$$
(this is easily proved for $z \in \dom{\opA^*}$ and thus remains true for $z \in H$ by density and continuity of $\Phi_T^*$ for any $T>0$).
Therefore, thanks to the assumption \eqref{dim ker PhiTmin finie} we know that
\begin{equation}\label{dim ker PhiTmin finie bis}
\dim \ker \Phi_{T}^* <+\infty, \quad \forall T \geq T_0.
\end{equation}
From now on, $T$ is fixed such that $T>T_0$.
Let $z \in \ker \Phi_T^*$.
To prove that $z \in \dom{A^*}$, we have to show that, for any sequence $t_n>0$ with $t_n \to 0$ as $n \to +\infty$, the sequence
$$u_n=\frac{S_{\opA}(t_n)^*z-z}{t_n}$$
converges in $H$ as $n \to +\infty$.
Let $\varepsilon \in (0, T-\Tmin]$ be fixed so that $T-\varepsilon \geq \Tmin$ and let $N \in \N$ be large enough so that $t_n<\varepsilon$ for every $n \geq N$.
Let us first show that
\begin{equation}\label{un dans kerPhi}
u_n \in \ker \Phi_{T-\varepsilon}^*, \quad \forall n \geq N.
\end{equation}
To this end, observe that, for $n \geq N$, we have (arguing by density as before)
$$
\int_0^{T-\varepsilon} \norm{\Phi_{T-\varepsilon}^*S_{\opA}(t_n)^*z(t)}_{U}^2 \, dt
=\int_{\varepsilon-t_n}^{T-t_n} \norm{\Phi_T^*z(t)}_{U}^2 \, dt.
$$
This shows that $S_{\opA}(t_n)^*z \in \ker \Phi_{T-\varepsilon}^*$ for $n \geq N$ since $z \in \ker \Phi_T^*$.
Thus, we have \eqref{un dans kerPhi}.
Let now $\mu \in \rho(\opA^*)$ be fixed and let us introduce the following norm on $\ker \Phi_{T-\varepsilon}^*$:
$$\norm{z}_{-1}=\norm{(\mu-\opA^*)^{-1}z}_H.$$
Since $(\mu-\opA^*)^{-1}z \in \dom{\opA^*}$, we have
$$
(\mu-\opA^*)^{-1}u_n
=\frac{S_{\opA}(t_n)^*-\Id}{t_n}(\mu-\opA^*)^{-1}z \xrightarrow[n \to +\infty]{}
\opA^*(\mu-\opA^*)^{-1}z \quad \mbox{ in } H.
$$
Therefore, $(u_n)_{n \geq N}$ is a Cauchy sequence in $\ker \Phi_{T-\varepsilon}^*$ for the norm $\norm{\cdot}_{-1}$.
Since $\ker \Phi_{T-\varepsilon}^*$ is finite dimensional (by \eqref{dim ker PhiTmin finie bis}), all the norms are equivalent on $\ker \Phi_{T-\varepsilon}^*$.
Thus, $(u_n)_{n \geq N}$ is then a Cauchy sequence for the usual norm $\norm{\cdot}_H$ as well and, as a result, converges for this norm.
This shows that $z \in \dom{\opA^*}$ and establishes \eqref{in the dom}.

Thanks to \eqref{in the dom} it is now easy to see that $\ker \Phi_T^*$ is stable by $\opA^*$ and that we have
\begin{equation}\label{ker PhiT dans ker Bs}
\ker \Phi_T^* \subset \ker \opB^*.
\end{equation}
Indeed, let $z \in \ker \Phi_T^*$, that is
\begin{equation}\label{z dans ker PhiT}
\Phi_T^*z(t)=0, \quad \aew t \in (0,T).
\end{equation}
On the one hand, since $z \in \dom{A^*}$, we can differentiate this identity to obtain (see e.g. \cite[Proposition 4.3.4]{TW09})
$$\Phi_T^*\opA^*z(t)=0, \quad \aew t \in (0,T),$$
that is $A^*z \in \ker \Phi_T^*$.
On the other hand, since $z \in \dom{A^*}$, the map $\Phi_T^*z:t \in [0,T] \mapsto \opB^*S_{\opA}(T-t)^*z \in U$ is continuous and we can take $t=T$ in \eqref{z dans ker PhiT} to obtain that $\opB^*z=\Phi_T^*z(T)=0$.

Consequently, the restriction $M$ of ${\opA}^*$ to $\ker \Phi_T^*$ is a linear operator from the finite dimensional space $\ker \Phi_T^*$ into itself.
Assume that $\ker \Phi_T^*\neq\ens{0}$ (otherwise the conclusion \eqref{ker PhiT} is clear from \eqref{ker PhiT supset}).
Therefore, $M$ is triangularizable in $\ker \Phi_T^*$ (here we use that $H$ is a complex Hilbert space).
In other words, $\ker \Phi_T^*$ is the direct sum of the root subspaces of $M$: for every $\lambda \in \spec{M}$, there exists $m(\lambda) \in \N^*$ such that
$$\ker \Phi_T^*=\bigoplus_{\lambda \in \spec{M}} \ker \left(\lambda-M\right)^{m(\lambda)}.$$
Finally, thanks to \eqref{ker PhiT dans ker Bs} we have $\spec{M} \subset \sigma_F$ and $\ker \left(\lambda-M\right)^{m(\lambda)} \subset E_\lambda$ for every $\lambda \in \spec{M}$.
\end{proof}

\begin{remark}
From the proof of Lemma \ref{lem fatt} we easily see that the equality \eqref{ker PhiT} remains valid for $T=\Tmin$ too if $\ker \Phi_{\Tmin}^* \subset \dom{A^*}$ (in addition to \eqref{dim ker PhiTmin finie}).
Let us also mention an alternative proof in the case $\sigma_F=\emptyset$.
Indeed, in this case, the same reasoning as at the end of the proof of Lemma \ref{lem fatt} easily shows that $\ker \Phi_{T_0}^* \cap \dom{{(A^*)}^{\infty}}=\ens{0}$.
By \cite[Proposition 6.1.9]{TW09} this implies that $\ker \Phi_T^*=\ens{0}$ for every $T>T_0$.
\end{remark}

Let us now conclude this section with the proof of Theorem \ref{thm pert 2}.

\begin{proof}[Proof of Theorem \ref{thm pert 2}]
The assumptions of Lemma \ref{lemme peetre} are satisfied for $L=\Phi_T^*$ and $P=\Psi_T^*-\Phi_T^*$ for every $T \geq \Tmin$.
Therefore, for every $T \geq \Tmin$, we have
$$\dim \ker \Phi_T^* <+\infty, \quad \left(\ker \Phi_T^*\right)^\perp=\range \Phi_T.$$
Applying now Lemma \ref{lem fatt} we obtain the desired conclusion.
\end{proof}

\section{An example}\label{sect 3}

Our results, especially Corollary \ref{thm pert fatt}, potentially have a lot of applications.
For instance, it can be used to investigate the exact controllability of some perturbed wave equations (by terms of order zero, as in \cite{Z91}), plate equations (by terms of order up to one, as in \cite{CT10}), Schr\"odinger equations, Korteweg-de Vries equations, transport equations (by non local spatial terms, as in \cite{CHO16}), etc.
In this section we focus on a particular example.

Let $\Omega \subset \R^N$ ($N \geq 1$) be an open bounded connected subset with boundary $\partial\Omega$ of class $C^2$ and let $\omega \subset \Omega$ be a non empty open subset.
Let $T>0$.
We consider the following Schr\"odinger equation with non local term:
\begin{equation}\label{schro nonloc}
\left\{\begin{array}{ll}
\ds y_t=i\Delta y +\int_\Omega k(\xi)y(t,\xi) \, d\xi +\indic_\omega(x) u(t,x)  &\mbox{ in } (0,T) \times \Omega, \\
\ds y =0 &\mbox{ on } (0,T) \times \partial\Omega, \\
\ds y(0)=y^0 &\mbox{ in } \Omega.
\end{array}\right.
\end{equation}
In \eqref{schro nonloc}, $y^0$ is the initial data, $y$ is the state and $u$ is the control.
$\indic_\omega$ denotes the function that is equal to $1$ in $\omega$ and $0$ outside.
The kernel $k$ is only assumed to be in $L^2(\Omega)$.
All the functions we consider in this section are complex-valued functions.

Let us recast \eqref{schro nonloc} as an abstract evolution system \eqref{syst abst}.
The state space $H$ and the control space $U$ are
$$H=U=L^2(\Omega).$$
The operator $\opA_{\K}:\dom{\opA_{\K}} \subset H \longrightarrow H$ is
$$
\opA_{\K} y=i\Delta y+\int_\Omega k(\xi)y(\xi) \, d\xi,
\qquad
\dom{\opA_{\K}}=H^2(\Omega) \cap H^1_0(\Omega),
$$
and the control operator $\opB:U \longrightarrow H$ is
$$\opB u=\indic_{\omega}u.$$
Clearly, $A_K$ splits up into $\opA_{\K}=\opA_{0}+\K$, where $\opA_{0}:\dom{\opA_{0}} \subset H \longrightarrow H$ is given by
$$
\opA_{0} y=i \Delta y,
\qquad
\dom{\opA_{0}}=\dom{\opA_{\K}},
$$
and $\K: H \longrightarrow H$ is given by
$$\K y=\int_\Omega k(\xi)y(\xi) \, d\xi.$$
The operator $\opA_0$ is skew-adjoint and therefore, by Stone's theorem, it is the generator of a unitary $C_0$-group on $H$.
On the other hand, it is clear that $\K$ is compact.
Finally, observe that $\opB$ is bounded and thus admissible.
Therefore, the assumptions \ref{hyp1 pAC} and \ref{hyp2 pAC} of Theorem \ref{thm pert AC} are satisfied.

We can now state the following simple (but new) consequence of Corollary \ref{thm pert fatt}:

\begin{theorem}\label{thm schro}
Assume that the Schr\"odinger equation $(A_0,B)$ is exactly controllable in time $\Tmin>0$.
Let us introduce the subset $\kappa \subset L^2(\Omega)$ defined by
$$\kappa=\ens{i\Delta \phi+\lambda \phi \quad \middle| \quad \lambda \in \C, \, \, \phi \in H^2(\Omega) \cap H^1_0(\Omega), \, \, \phi=0 \mbox{ in } \omega, \, \, \int_\Omega \phi(\xi) \, d\xi=1}.$$
If $k \not\in \kappa$, then, for every $T>\Tmin$, \eqref{schro nonloc} is exactly controllable in time $T$.
Conversely, if $k \in \kappa$, then, for every $T>0$, \eqref{schro nonloc} is not approximately controllable in time $T$.
\end{theorem}

Let us point out that $\kappa$ is never empty, unless $\omega=\Omega$.
For explicit conditions of the exact controllability of the Schr\"odinger equation $(A_0,B)$, we refer to the survey \cite{L14} and the references therein.
It is for instance nowadays well-known that this equation is exactly controllable in time $T_0$ for every $T_0>0$ if the control domain $\omega$ satisfies the so-called Geometric Control Condition.

\begin{proof}[Proof of Theorem \ref{thm schro}]
Assume first that $k \not\in \kappa$.
We are going to apply Corollary \ref{thm pert fatt}.
We only have to check the Fattorini-Hautus test corresponding to \eqref{schro nonloc}.
It is clear that $B^*=B$ and $A_K^*=-A_0+K^*$ with $\dom{A_K^*}=\dom{A_K}$ and
$$K^*z=k \int_{\Omega} z(\xi) \, d\xi.$$
Let then $\lambda \in \C$ and $z \in H^2(\Omega) \cap H^1_0(\Omega)$ be such that
\begin{equation}\label{fatt test schro nonloc}
\left\{\begin{array}{ll}
\ds -i\Delta z+k \int _\Omega z(\xi) \, d\xi=\lambda z &\mbox{ in } \Omega, \\
\ds z=0 &\mbox{ in } \omega,
\end{array}\right.
\end{equation}
and let us show that this implies that $z=0$ in $\Omega$.
Since $k \not\in \kappa$, by the very definition of $\kappa$, we obtain that the constant $\int _\Omega z(\xi) \, d\xi$ is necessarily equal to zero.
Thus, $z \in H^2(\Omega) \cap H^1_0(\Omega)$ satisfies
$$
\left\{\begin{array}{ll}
\ds -i\Delta z=\lambda z &\mbox{ in } \Omega, \\
\ds z=0 &\mbox{ in } \omega.
\end{array}\right.
$$
As it is well-known, this implies that $z=0$ in $\Omega$ (we recall that $\Omega$ is connected).

Conversely, it is easy to see that, if $k \in \kappa$, then there exists a non zero solution $z$ to \eqref{fatt test schro nonloc}.
\end{proof}

\begin{remark}
Obviously, we can consider more general kernels $k \in L^2(\Omega \times \Omega)$ that depend on the $x$ variable as well.
However, the condition provided by the Fattorini-Hautus test may be less explicit than the one presented in Theorem \ref{thm schro}.
We also chose to apply Corollary \ref{thm pert fatt} instead of Theorem \ref{thm pert AC} because it does not seem so easy to check the approximate controllability assumption \ref{item ac}.
One case where it is not difficult to directly check this condition is when the kernel satisfies the analyticity assumption (3) of \cite{FCLZ16}, namely that, for every $y \in L^2(\Omega)$,
\begin{equation}\label{hyp FCLZ}
x \in \Omega \longmapsto \int_{\Omega} k(\xi,x)y(\xi) \, d\xi \mbox{ is analytic.}
\end{equation}
Note that \eqref{hyp FCLZ} is not satisfied here since we only assumed the kernel to be $L^2$.
\end{remark}

\section{Additional comments}\label{sect 4}

Let us conclude this article with some final important remarks.

Our main results were obtained as a combination of three ingredients, namely:
\begin{enumerate}[1)]
\item\label{step 1}
the desired observability inequality with a ``compact error'' \eqref{obs estim cpct err} (obtained in this article as a consequence of the assumptions \ref{hyp1 pAC} and \ref{hyp2 pAC}),
\item\label{step 2}
the Peetre lemma,
\item\label{step 3}
the general Lemma \ref{lem fatt}.
\end{enumerate}

In the course of this paper we have seen that this procedure is general and we have actually establish the following abstract result:

\begin{theorem}\label{main thm}
Let $H$ and $U$ be two complex Hilbert spaces.
Let $\opA:\dom{\opA} \subset H \longrightarrow H$ be the generator of a $C_0$-semigroup on $H$ and let $\opB \in \lin{U,\dom{\opA^*}'}$ be admissible for $\opA$.
For $T \geq 0$ let $\Phi_T \in \lin{L^2(0,+\infty;U),H}$ be the input map of $(A,B)$.
Assume that there exist a complex Banach space $V$ and a compact operator $P \in \lin{H,V}$ such that, for some $T_0>0$ and $\alpha>0$, we have
\begin{equation}\label{obs ineq plus cpct}
\alpha\norm{z}_H^2 \leq \int_0^{T_0} \norm{\Phi_{T_0}^*z(t)}_U^2 \, dt+\norm{P z}_{V}^2, \quad \forall z \in H.
\end{equation}
Then, the set $\sigma_F$ defined by \eqref{fatt spec} is finite, the subspace $E_\lambda$ defined by \eqref{Elambda} is finite dimensional for every $\lambda \in \C$, and we have
$$\range \Phi_T =\left(\bigoplus_{\lambda \in \sigma_F} E_\lambda\right)^{\perp}, \quad \forall T>\Tmin.$$
In particular, $(A,B)$ is exactly controllable in time $T$ for every $T>T_0$ if $(A,B)$ satisfies the Fattorini-Hautus test.
\end{theorem}

In Theorem \ref{thm pert 2}, we have provided some conditions on $(A,B)$ to fulfill the assumption \eqref{obs ineq plus cpct}.
However, there are many other techniques that can be used to establish \eqref{obs ineq plus cpct}.
For instance, an easy way to obtain such estimates is to use the so-called multiplier method and then to check to that the ``error term'' is indeed compact thanks to some compactness results such as the theorem of Aubin-Lions-Simon.
This is precisely what is achieved for the plate equation in \cite[Th\'eor\`eme IV.3.7]{L88} (there, the estimate (3.127) plays the role of \eqref{obs ineq plus cpct}) or for a Korteweg-de Vries in \cite[Section 3]{R97} (see the estimate (3.13)).
These are just some examples of equations where the multiplier method leads to \eqref{obs ineq plus cpct}.

In conclusion, Theorem \ref{main thm} is not only useful to investigate the controllability of compactly perturbed controlled systems, but it is also useful to establish the controllability of numerous unperturbed equations themselves.
This abstract result can be considered as a formalization of the compactness-uniqueness method used for the controllability of second order systems in \cite{L88} and used to study the stability of some hyperbolic equations in \cite{RT74}.

\begin{remark}
With the notation of Lemma \ref{lemme peetre}, we can also show that \ref{hyp a} is equivalent to the following condition:
\begin{enumerate}[(c)]
\item\label{hyp c}
There exist two subspaces $M,N \subset B_1$ such that
$$B_1=M \oplus N,$$
with, for some $c>0$,
\begin{equation}\label{obs ineq plus cpct 2}
c \norm{z}_{B_1} \leq \norm{Lz}_{B_2}, \quad \forall z \in M,
\end{equation}
and $\dim N<+\infty$ with the projection from $B_1$ onto $N$ being continuous.
\end{enumerate}
It is not difficult to see that \ref{hyp c} implies \ref{hyp a} with $B_3=B_1$ and $P$ being the projection from $B_1$ onto $N$, which is continuous by assumption and compact since its range $N$ is finite dimensional.
On the other hand, that \ref{hyp a} implies \ref{hyp c} is actually contained in the proof of Lemma \ref{lemme peetre} (see the proof of ii) $\Longrightarrow$ i) of \cite[Lemma 3]{P61}).

This provides an equivalent formulation of the assumption \eqref{obs ineq plus cpct} which can be useful for some problems.
Notably, we see that Theorem \ref{main thm} generalizes \cite[Theorem 5.2]{K94} by removing some unnecessary spectral assumptions made on the operator $A$.
It is also more general than \cite[Lemma 1.1]{LR17} since this lemma requires that the corresponding right-hand side in \eqref{obs ineq plus cpct 2} defines a norm on $B_1$, which is equivalent to assuming that $(A,B)$ is approximately controllable in time $T_0$ (and we recall that this property is stronger than the Fattorini-Hautus test).
\end{remark}

\appendix
\section{Proof of the estimate \eqref{key estim}}\label{appendix}

This appendix is devoted to a proof of the estimate \eqref{key estim} that is used in the proof of Lemma \ref{lem cpct}.
It is largely inspired by \cite[Proposition 3.3]{H05}.

Let us recall our framework.
$H$ and $U$ are two Hilbert spaces.
$\opA:\dom{\opA} \subset H \longrightarrow H$ is the generator of a $C_0$-semigroup $(S_{\opA}(t))_{t \geq 0}$ on $H$ and $\opB \in \lin{U,\dom{\opA^*}'}$ is admissible for $\opA$.
Let us recall the following dual characterization of admissibility: $B$ is admissible for $A$ if, and only if, for some (and hence all) $T>0$, there exists $\gamma>0$ such that
\begin{equation}\label{B admiss}
\int_0^T \norm{B^*S_A(T-t)^*z}_U^2 \, dt \leq \gamma \norm{z}_H^2, \quad \forall z \in \dom{A^*}.
\end{equation}
Let us now introduce for $n \in \N$ large enough ($n>\omega_0$, where $\omega_0 \in \R$ is the growth bound of $A$) the Yosida-like approximations $C_n \in \lin{H,U}$ defined by
$$C_n z=nB^*(n-A^*)^{-1}z, \quad \forall z \in H.$$
Let us recall that (see e.g. \cite[Lemma II.3.4]{EN})
\begin{equation}\label{conv reso}
n(n-A^*)^{-1}z \xrightarrow[n \to +\infty]{} z \quad \mbox{ in } H, \quad \forall z \in H.
\end{equation}
This implies in particular that
\begin{equation}\label{Cn conv}
C_n z \xrightarrow[n \to +\infty]{} B^*z \quad \mbox{ in } U, \quad \forall z \in \dom{A^*},
\end{equation}
since for every $z \in \dom{A^*}$ we have
\begin{multline*}
\norm{C_n z-B^*z}_U \\
\leq \norm{B^*}_{\lin{\dom{A^*},U}} \left(\norm{n(n-A^*)^{-1}A^*z-A^*z}_{H}+\norm{n(n-A^*)^{-1}z-z}_{H} \right).
\end{multline*}
Let $T>0$ be fixed.
For $f \in L^2(0,T;H)$, let us denote by $S_A^* \ast f \in L^2(0,T;H)$ the function defined for every $t \in (0,T)$ by
$$(S_A^* \ast f)(t)=\int_0^t S_{\opA}(t-s)^*f(s) \,ds.$$
Using the Cauchy-Schwarz inequality we have
$$
\int_0^T \norm{C_n(S_A^* \ast f)(t)}_U^2 \, dt
\leq T \int_0^T \int_0^t  \norm{B^*S_{\opA}(t-s)^*n(n-A^*)^{-1}f(s)}_U^2  \,ds  \, dt.
$$
Using Fubini's theorem we obtain
$$
\int_0^T \norm{C_n(S_A^* \ast f)(t)}_U^2 \, dt
\leq T \int_0^T \int_0^T  \norm{B^*S_{\opA}(t)^*n(n-A^*)^{-1}f(s)}_U^2  \,dt  \, ds.
$$
Using now the admissibility of $B$ (see \eqref{B admiss}) this gives
\begin{equation}\label{key estim Cn}
\int_0^T \norm{C_n(S_A^* \ast f)(t)}_U^2 \, dt \leq T\gamma \int_0^T \norm{n(n-A^*)^{-1}f(s)}_H^2  \,ds.
\end{equation}
Let us now remark that $S_A^* \ast f \in L^2(0,T;\dom{A^*})$ for $f \in C^1([0,T];H)$.
Using then \eqref{Cn conv} and \eqref{conv reso}, we see that the Lebesgue's dominated convergence theorem applies and that we can pass to the limit $n \to +\infty$ in \eqref{key estim Cn} to finally obtain the desired estimate.

\bibliographystyle{siam}
\bibliography{biblio}

\medskip
\medskip

\end{document}